\documentclass{gtart}

\def\ifplaintex{\expandafter\ifx\csname documentclass\endcsname\relax}

\def\gtp{{\mathsurround=0pt\it $\cal G\mskip-2mu$eometry \&\ 
$\cal T\!\!$opology $\cal P\!$ublications}}  

\def\Addressesr{\bigskip
{\small \parskip 0pt \leftskip 0pt \rightskip 0pt plus 1fil \def\\{\par}
\sl\theaddress\par
\medskip
\rm Email:\stdspace\tt\theemail\hfill\rm Received:\qua\receiveddate \par}}

\def\recd{{\small Received:\qua\receiveddate\ifx\reviseddate\relax
\else\qquad Revised:\qua\reviseddate\fi\par}} 

\def\lognumber#1{\def\thelognumber{#1}}
\def\volumenumber#1{\def\thevolumenumber{#1}}
\def\volumeyear#1{\def\thevolumeyear{#1}}
\def\papernumber#1{\def\thepapernumber{#1}}
\def\pagenumbers#1#2{\def\startpage{#1}\def\finishpage{#2}}
\def\published#1{\def\publishdate{#1}}

\def\received#1{\def\receiveddate{#1}}

\def\accepted#1{\def\accepteddate{#1}}

\def\asciiaddress#1{\def\theasciiaddress{#1}}

\long\def\asciiabstract#1{\long\def\theasciiabstract{#1}}

\let\\\par\let\thelognumber\relax\let\thevolumenumber\relax
\let\thepapernumber\relax\let\thevolumeyear\relax\let\startpage\relax
\let\finishpage\relax\let\publishdate\relax\let\receiveddate\relax
\let\reviseddate\relax\let\accepteddate\relax\let\theasciititle\relax
\let\theasciiauthors\relax\let\theasciiaddress\relax
\let\theasciiabstract\relax

\let\theasciiemail\relax

\ifplaintex
\font\logobig=cmssbx10 scaled 3836
\font\logomed=cmssbx10 scaled 2557
\else
\font\logobig=cmssbx10 scaled 4200
\font\logomed=cmssbx10 scaled 2800
\fi

\long\def\makeagttitle{   
\count0=\startpage
\agt\hfill      
\hbox to 45truept{\vbox to 0pt{\vglue -13truept{\logomed A\kern -.37em{\logobig 
T}\kern -.38em G}\vss}\hss}
\break
{\small Volume \thevolumenumber\ (\thevolumeyear)
\startpage--\finishpage\nl
Published: \publishdate}

\vglue .25truein

{\parskip=0pt\leftskip 0pt plus
1fil\def\\{\par\smallskip}{\Large\bf\thetitle}\par\medskip} \vglue
0.05truein

{\parskip=0pt\leftskip 0pt plus 1fil\def\\{\par}{\sc\theauthors}
\par\medskip}%
 
\vglue 0.03truein 

{\small\leftskip 25truept\rightskip 25truept{\bf Abstract}\stdspace\theabstract

{\bf AMS Classification}\stdspace\theprimaryclass
\ifx\thesecondaryclass\relax\else; \thesecondaryclass\fi\par
{\bf Keywords}\stdspace \thekeywords\par}\vglue 7truept

}   

\ifplaintex
\font\phead=cmsl9 scaled 950
\font\pnum=cmbx10 scaled 913
\font\pfoot=cmsl9 scaled 950
\headline{\vbox to 0pt{\vskip -4.5mm\line{\small\phead\ifnum
\count0=\startpage ISSN 1472-2739 (on-line) 1472-2747 (printed)
\hfill {\pnum\folio}\else\ifodd\count0\def\\{ }%
\ifx\theshorttitle\relax\thetitle\else\theshorttitle\fi\hfill{\pnum\folio}
\else\def\\{ and }{\pnum\folio}\hfill\ifx\theshortauthors\relax\theauthors
\else\theshortauthors\fi\fi\fi}\vss}}
\footline{\vbox to 0pt{\vglue 0mm\line{\small\pfoot\ifnum\count0=\startpage
\copyright\ \gtp\hfill\else
\agt, Volume \thevolumenumber\ (\thevolumeyear)\hfill\fi}\vss}}
\else
\font=cmsl9 scaled 1050
\font\lnum=cmbx10 
\font=cmsl9 scaled 1050
\makeatletter
\def\@oddhead{{\small\lhead\ifnum\count0=\startpage ISSN 1472-2739 
(on-line) 1472-2747 (printed)\hfill {\lnum\number\count0}\else\ifodd\count0
\def\\{ }\ifx\theshorttitle\relax \thetitle \else\theshorttitle\fi\hfill
{\lnum\number\count0}\else\def\\{ and }{\lnum\number\count0}
\hfill\ifx\theshortauthors\relax 
\theauthors\else\theshortauthors\fi\fi\fi}}\def\@evenhead{\@oddhead}
\def\@oddfoot{\small\lfoot\ifnum\count0=\startpage\copyright\ \gtp\hfill\else
\agt, Volume \thevolumenumber\ (\thevolumeyear)\hfill\fi}
\def\@evenfoot{\@oddfoot}
\makeatother
\fi
\let\maketitlepage\makeagttitle
\let\makeshorttitle\maketitlepage
\let\maketitle\maketitlepage

\newwrite\gtoutfile
\long\gdef\makeheadfile{  
{\def\\{, }\def\s{ }
\immediate\openout\gtoutfile head.xxx
\immediate\write\gtoutfile{To: math@arxiv.org}
\immediate\write\gtoutfile{Subject: put OR rep NNNNN:ppppp}
\immediate\write\gtoutfile{--text follows this line--}
\immediate\write\gtoutfile{Proxy-for: \ifx\theasciiauthors\relax
\theauthors\else\theasciiauthors\fi\s<\ifx\theasciiemail\relax\theemail\else\theasciiemail\fi>}
\immediate\write\gtoutfile{\noexpand\\}
\immediate\write\gtoutfile{Authors: \ifx\theasciiauthors\relax
\theauthors\else\theasciiauthors\fi}
{\def\\{ }\immediate\write\gtoutfile{Title: \ifx\theasciititle\relax
\thetitle\else\theasciititle\fi}}
\immediate\write\gtoutfile{Subj-class: GT or SG, GR etc}
\immediate\write\gtoutfile{MSC-class: \theprimaryclass\ifx\thesecondaryclass\relax\else, \thesecondaryclass\fi}
\immediate\write\gtoutfile{Journal-ref: Algebr. Geom. Topol. \thevolumenumber\s
(\thevolumeyear) \startpage-\finishpage}
\immediate\write\gtoutfile{Comments: Published by Algebraic and
Geometric Topology at}
\immediate\write\gtoutfile{\s\s\s  http://www.maths.warwick.ac.uk/agt/AGTVol\thevolumenumber/agt-\thevolumenumber-\thepapernumber.abs.html}
\immediate\write\gtoutfile{\noexpand\\}
\immediate\write\gtoutfile{}
\ifx\theasciiabstract\relax
\immediate\write\gtoutfile{\theabstract}\else
\immediate\write\gtoutfile{\theasciiabstract}\fi
\immediate\write\gtoutfile{}
\immediate\write\gtoutfile{\noexpand\\}
\immediate\write\gtoutfile{}
\immediate\closeout\gtoutfile}}  

\def\maketitlepage{\makeagttitle\makeheadfile}
\let\makeshorttitle\maketitlepage
\let\maketitle\maketitlepage

\def\ifplaintex{\expandafter\ifx\csname documentclass\endcsname\relax}

\def\gtp{{\mathsurround=0pt\it $\cal G\mskip-2mu$eometry \&\ 
$\cal T\!\!$opology $\cal P\!$ublications}}  

\def\Addressesr{\bigskip
{\small \parskip 0pt \leftskip 0pt \rightskip 0pt plus 1fil \def\\{\par}
\sl\theaddress\par
\medskip
\rm Email:\stdspace\tt\theemail\hfill\rm Received:\qua\receiveddate \par}}

\def\recd{{\small Received:\qua\receiveddate\ifx\reviseddate\relax
\else\qquad Revised:\qua\reviseddate\fi\par}} 

\def\lognumber#1{\def\thelognumber{#1}}
\def\volumenumber#1{\def\thevolumenumber{#1}}
\def\volumeyear#1{\def\thevolumeyear{#1}}
\def\papernumber#1{\def\thepapernumber{#1}}
\def\pagenumbers#1#2{\def\startpage{#1}\def\finishpage{#2}}
\def\published#1{\def\publishdate{#1}}

\def\received#1{\def\receiveddate{#1}}

\def\accepted#1{\def\accepteddate{#1}}

\def\asciiaddress#1{\def\theasciiaddress{#1}}

\long\def\asciiabstract#1{\long\def\theasciiabstract{#1}}

\let\\\par\let\thelognumber\relax\let\thevolumenumber\relax
\let\thepapernumber\relax\let\thevolumeyear\relax\let\startpage\relax
\let\finishpage\relax\let\publishdate\relax\let\receiveddate\relax
\let\reviseddate\relax\let\accepteddate\relax\let\theasciititle\relax
\let\theasciiauthors\relax\let\theasciiaddress\relax
\let\theasciiabstract\relax

\let\theasciiemail\relax

\ifplaintex
\font\logobig=cmssbx10 scaled 3836
\font\logomed=cmssbx10 scaled 2557
\else
\font\logobig=cmssbx10 scaled 4200
\font\logomed=cmssbx10 scaled 2800
\fi

\long\def\makeagttitle{   
\count0=\startpage
\agt\hfill      
\hbox to 45truept{\vbox to 0pt{\vglue -13truept{\logomed A\kern -.37em{\logobig 
T}\kern -.38em G}\vss}\hss}
\break
{\small Volume \thevolumenumber\ (\thevolumeyear)
\startpage--\finishpage\nl
Published: \publishdate}

\vglue .25truein

{\parskip=0pt\leftskip 0pt plus
1fil\def\\{\par\smallskip}{\Large\bf\thetitle}\par\medskip} \vglue
0.05truein

{\parskip=0pt\leftskip 0pt plus 1fil\def\\{\par}{\sc\theauthors}
\par\medskip}%
 
\vglue 0.03truein 

{\small\leftskip 25truept\rightskip 25truept{\bf Abstract}\stdspace\theabstract

{\bf AMS Classification}\stdspace\theprimaryclass
\ifx\thesecondaryclass\relax\else; \thesecondaryclass\fi\par
{\bf Keywords}\stdspace \thekeywords\par}\vglue 7truept

}   

\ifplaintex
\font\phead=cmsl9 scaled 950
\font\pnum=cmbx10 scaled 913
\font\pfoot=cmsl9 scaled 950
\headline{\vbox to 0pt{\vskip -4.5mm\line{\small\phead\ifnum
\count0=\startpage ISSN 1472-2739 (on-line) 1472-2747 (printed)
\hfill {\pnum\folio}\else\ifodd\count0\def\\{ }%
\ifx\theshorttitle\relax\thetitle\else\theshorttitle\fi\hfill{\pnum\folio}
\else\def\\{ and }{\pnum\folio}\hfill\ifx\theshortauthors\relax\theauthors
\else\theshortauthors\fi\fi\fi}\vss}}
\footline{\vbox to 0pt{\vglue 0mm\line{\small\pfoot\ifnum\count0=\startpage
\copyright\ \gtp\hfill\else
\agt, Volume \thevolumenumber\ (\thevolumeyear)\hfill\fi}\vss}}
\else
\font=cmsl9 scaled 1050
\font\lnum=cmbx10 
\font=cmsl9 scaled 1050
\makeatletter
\def\@oddhead{{\small\lhead\ifnum\count0=\startpage ISSN 1472-2739 
(on-line) 1472-2747 (printed)\hfill {\lnum\number\count0}\else\ifodd\count0
\def\\{ }\ifx\theshorttitle\relax \thetitle \else\theshorttitle\fi\hfill
{\lnum\number\count0}\else\def\\{ and }{\lnum\number\count0}
\hfill\ifx\theshortauthors\relax 
\theauthors\else\theshortauthors\fi\fi\fi}}\def\@evenhead{\@oddhead}
\def\@oddfoot{\small\lfoot\ifnum\count0=\startpage\copyright\ \gtp\hfill\else
\agt, Volume \thevolumenumber\ (\thevolumeyear)\hfill\fi}
\def\@evenfoot{\@oddfoot}
\makeatother
\fi
\let\maketitlepage\makeagttitle
\let\makeshorttitle\maketitlepage
\let\maketitle\maketitlepage

\newwrite\gtoutfile
\long\gdef\makeheadfile{  
{\def\\{, }\def\s{ }
\immediate\openout\gtoutfile head.xxx
\immediate\write\gtoutfile{To: math@arxiv.org}
\immediate\write\gtoutfile{Subject: put OR rep NNNNN:ppppp}
\immediate\write\gtoutfile{--text follows this line--}
\immediate\write\gtoutfile{Proxy-for: \ifx\theasciiauthors\relax
\theauthors\else\theasciiauthors\fi\s<\ifx\theasciiemail\relax\theemail\else\theasciiemail\fi>}
\immediate\write\gtoutfile{\noexpand\\}
\immediate\write\gtoutfile{Authors: \ifx\theasciiauthors\relax
\theauthors\else\theasciiauthors\fi}
{\def\\{ }\immediate\write\gtoutfile{Title: \ifx\theasciititle\relax
\thetitle\else\theasciititle\fi}}
\immediate\write\gtoutfile{Subj-class: GT or SG, GR etc}
\immediate\write\gtoutfile{MSC-class: \theprimaryclass\ifx\thesecondaryclass\relax\else, \thesecondaryclass\fi}
\immediate\write\gtoutfile{Journal-ref: Algebr. Geom. Topol. \thevolumenumber\s
(\thevolumeyear) \startpage-\finishpage}
\immediate\write\gtoutfile{Comments: Published by Algebraic and
Geometric Topology at}
\immediate\write\gtoutfile{\s\s\s  http://www.maths.warwick.ac.uk/agt/AGTVol\thevolumenumber/agt-\thevolumenumber-\thepapernumber.abs.html}
\immediate\write\gtoutfile{\noexpand\\}
\immediate\write\gtoutfile{}
\ifx\theasciiabstract\relax
\immediate\write\gtoutfile{\theabstract}\else
\immediate\write\gtoutfile{\theasciiabstract}\fi
\immediate\write\gtoutfile{}
\immediate\write\gtoutfile{\noexpand\\}
\immediate\write\gtoutfile{}
\immediate\closeout\gtoutfile}}  

\def\maketitlepage{\makeagttitle\makeheadfile}
\let\makeshorttitle\maketitlepage
\let\maketitle\maketitlepage

\def\ifplaintex{\expandafter\ifx\csname documentclass\endcsname\relax}

\def\gtp{{\mathsurround=0pt\it $\cal G\mskip-2mu$eometry \&\ 
$\cal T\!\!$opology $\cal P\!$ublications}}  

\def\Addressesr{\bigskip
{\small \parskip 0pt \leftskip 0pt \rightskip 0pt plus 1fil \def\\{\par}
\sl\theaddress\par
\medskip
\rm Email:\stdspace\tt\theemail\hfill\rm Received:\qua\receiveddate \par}}

\def\recd{{\small Received:\qua\receiveddate\ifx\reviseddate\relax
\else\qquad Revised:\qua\reviseddate\fi\par}} 

\def\lognumber#1{\def\thelognumber{#1}}
\def\volumenumber#1{\def\thevolumenumber{#1}}
\def\volumeyear#1{\def\thevolumeyear{#1}}
\def\papernumber#1{\def\thepapernumber{#1}}
\def\pagenumbers#1#2{\def\startpage{#1}\def\finishpage{#2}}
\def\published#1{\def\publishdate{#1}}

\def\received#1{\def\receiveddate{#1}}

\def\accepted#1{\def\accepteddate{#1}}

\def\asciiaddress#1{\def\theasciiaddress{#1}}

\long\def\asciiabstract#1{\long\def\theasciiabstract{#1}}

\let\\\par\let\thelognumber\relax\let\thevolumenumber\relax
\let\thepapernumber\relax\let\thevolumeyear\relax\let\startpage\relax
\let\finishpage\relax\let\publishdate\relax\let\receiveddate\relax
\let\reviseddate\relax\let\accepteddate\relax\let\theasciititle\relax
\let\theasciiauthors\relax\let\theasciiaddress\relax
\let\theasciiabstract\relax

\let\theasciiemail\relax

\ifplaintex
\font\logobig=cmssbx10 scaled 3836
\font\logomed=cmssbx10 scaled 2557
\else
\font\logobig=cmssbx10 scaled 4200
\font\logomed=cmssbx10 scaled 2800
\fi

\long\def\makeagttitle{   
\count0=\startpage
\agt\hfill      
\hbox to 45truept{\vbox to 0pt{\vglue -13truept{\logomed A\kern -.37em{\logobig 
T}\kern -.38em G}\vss}\hss}
\break
{\small Volume \thevolumenumber\ (\thevolumeyear)
\startpage--\finishpage\nl
Published: \publishdate}

\vglue .25truein

{\parskip=0pt\leftskip 0pt plus
1fil\def\\{\par\smallskip}{\Large\bf\thetitle}\par\medskip} \vglue
0.05truein

{\parskip=0pt\leftskip 0pt plus 1fil\def\\{\par}{\sc\theauthors}
\par\medskip}%
 
\vglue 0.03truein 

{\small\leftskip 25truept\rightskip 25truept{\bf Abstract}\stdspace\theabstract

{\bf AMS Classification}\stdspace\theprimaryclass
\ifx\thesecondaryclass\relax\else; \thesecondaryclass\fi\par
{\bf Keywords}\stdspace \thekeywords\par}\vglue 7truept

}   

\ifplaintex
\font\phead=cmsl9 scaled 950
\font\pnum=cmbx10 scaled 913
\font\pfoot=cmsl9 scaled 950
\headline{\vbox to 0pt{\vskip -4.5mm\line{\small\phead\ifnum
\count0=\startpage ISSN 1472-2739 (on-line) 1472-2747 (printed)
\hfill {\pnum\folio}\else\ifodd\count0\def\\{ }%
\ifx\theshorttitle\relax\thetitle\else\theshorttitle\fi\hfill{\pnum\folio}
\else\def\\{ and }{\pnum\folio}\hfill\ifx\theshortauthors\relax\theauthors
\else\theshortauthors\fi\fi\fi}\vss}}
\footline{\vbox to 0pt{\vglue 0mm\line{\small\pfoot\ifnum\count0=\startpage
\copyright\ \gtp\hfill\else
\agt, Volume \thevolumenumber\ (\thevolumeyear)\hfill\fi}\vss}}
\else
\font=cmsl9 scaled 1050
\font\lnum=cmbx10 
\font=cmsl9 scaled 1050
\makeatletter
\def\@oddhead{{\small\lhead\ifnum\count0=\startpage ISSN 1472-2739 
(on-line) 1472-2747 (printed)\hfill {\lnum\number\count0}\else\ifodd\count0
\def\\{ }\ifx\theshorttitle\relax \thetitle \else\theshorttitle\fi\hfill
{\lnum\number\count0}\else\def\\{ and }{\lnum\number\count0}
\hfill\ifx\theshortauthors\relax 
\theauthors\else\theshortauthors\fi\fi\fi}}\def\@evenhead{\@oddhead}
\def\@oddfoot{\small\lfoot\ifnum\count0=\startpage\copyright\ \gtp\hfill\else
\agt, Volume \thevolumenumber\ (\thevolumeyear)\hfill\fi}
\def\@evenfoot{\@oddfoot}
\makeatother
\fi
\let\maketitlepage\makeagttitle
\let\makeshorttitle\maketitlepage
\let\maketitle\maketitlepage

\newwrite\gtoutfile
\long\gdef\makeheadfile{  
{\def\\{, }\def\s{ }
\immediate\openout\gtoutfile head.xxx
\immediate\write\gtoutfile{To: math@arxiv.org}
\immediate\write\gtoutfile{Subject: put OR rep NNNNN:ppppp}
\immediate\write\gtoutfile{--text follows this line--}
\immediate\write\gtoutfile{Proxy-for: \ifx\theasciiauthors\relax
\theauthors\else\theasciiauthors\fi\s<\ifx\theasciiemail\relax\theemail\else\theasciiemail\fi>}
\immediate\write\gtoutfile{\noexpand\\}
\immediate\write\gtoutfile{Authors: \ifx\theasciiauthors\relax
\theauthors\else\theasciiauthors\fi}
{\def\\{ }\immediate\write\gtoutfile{Title: \ifx\theasciititle\relax
\thetitle\else\theasciititle\fi}}
\immediate\write\gtoutfile{Subj-class: GT or SG, GR etc}
\immediate\write\gtoutfile{MSC-class: \theprimaryclass\ifx\thesecondaryclass\relax\else, \thesecondaryclass\fi}
\immediate\write\gtoutfile{Journal-ref: Algebr. Geom. Topol. \thevolumenumber\s
(\thevolumeyear) \startpage-\finishpage}
\immediate\write\gtoutfile{Comments: Published by Algebraic and
Geometric Topology at}
\immediate\write\gtoutfile{\s\s\s  http://www.maths.warwick.ac.uk/agt/AGTVol\thevolumenumber/agt-\thevolumenumber-\thepapernumber.abs.html}
\immediate\write\gtoutfile{\noexpand\\}
\immediate\write\gtoutfile{}
\ifx\theasciiabstract\relax
\immediate\write\gtoutfile{\theabstract}\else
\immediate\write\gtoutfile{\theasciiabstract}\fi
\immediate\write\gtoutfile{}
\immediate\write\gtoutfile{\noexpand\\}
\immediate\write\gtoutfile{}
\immediate\closeout\gtoutfile}}  

\def\maketitlepage{\makeagttitle\makeheadfile}
\let\makeshorttitle\maketitlepage
\let\maketitle\maketitlepage

\def\ifplaintex{\expandafter\ifx\csname documentclass\endcsname\relax}

\def\gtp{{\mathsurround=0pt\it $\cal G\mskip-2mu$eometry \&\ 
$\cal T\!\!$opology $\cal P\!$ublications}}  

\def\Addressesr{\bigskip
{\small \parskip 0pt \leftskip 0pt \rightskip 0pt plus 1fil \def\\{\par}
\sl\theaddress\par
\medskip
\rm Email:\stdspace\tt\theemail\hfill\rm Received:\qua\receiveddate \par}}

\def\recd{{\small Received:\qua\receiveddate\ifx\reviseddate\relax
\else\qquad Revised:\qua\reviseddate\fi\par}} 

\def\lognumber#1{\def\thelognumber{#1}}
\def\volumenumber#1{\def\thevolumenumber{#1}}
\def\volumeyear#1{\def\thevolumeyear{#1}}
\def\papernumber#1{\def\thepapernumber{#1}}
\def\pagenumbers#1#2{\def\startpage{#1}\def\finishpage{#2}}
\def\published#1{\def\publishdate{#1}}

\def\received#1{\def\receiveddate{#1}}

\def\accepted#1{\def\accepteddate{#1}}

\def\asciiaddress#1{\def\theasciiaddress{#1}}

\long\def\asciiabstract#1{\long\def\theasciiabstract{#1}}

\let\\\par\let\thelognumber\relax\let\thevolumenumber\relax
\let\thepapernumber\relax\let\thevolumeyear\relax\let\startpage\relax
\let\finishpage\relax\let\publishdate\relax\let\receiveddate\relax
\let\reviseddate\relax\let\accepteddate\relax\let\theasciititle\relax
\let\theasciiauthors\relax\let\theasciiaddress\relax
\let\theasciiabstract\relax

\let\theasciiemail\relax

\ifplaintex
\font\logobig=cmssbx10 scaled 3836
\font\logomed=cmssbx10 scaled 2557
\else
\font\logobig=cmssbx10 scaled 4200
\font\logomed=cmssbx10 scaled 2800
\fi

\long\def\makeagttitle{   
\count0=\startpage
\agt\hfill      
\hbox to 45truept{\vbox to 0pt{\vglue -13truept{\logomed A\kern -.37em{\logobig 
T}\kern -.38em G}\vss}\hss}
\break
{\small Volume \thevolumenumber\ (\thevolumeyear)
\startpage--\finishpage\nl
Published: \publishdate}

\vglue .25truein

{\parskip=0pt\leftskip 0pt plus
1fil\def\\{\par\smallskip}{\Large\bf\thetitle}\par\medskip} \vglue
0.05truein

{\parskip=0pt\leftskip 0pt plus 1fil\def\\{\par}{\sc\theauthors}
\par\medskip}%
 
\vglue 0.03truein 

{\small\leftskip 25truept\rightskip 25truept{\bf Abstract}\stdspace\theabstract

{\bf AMS Classification}\stdspace\theprimaryclass
\ifx\thesecondaryclass\relax\else; \thesecondaryclass\fi\par
{\bf Keywords}\stdspace \thekeywords\par}\vglue 7truept

}   

\ifplaintex
\font\phead=cmsl9 scaled 950
\font\pnum=cmbx10 scaled 913
\font\pfoot=cmsl9 scaled 950
\headline{\vbox to 0pt{\vskip -4.5mm\line{\small\phead\ifnum
\count0=\startpage ISSN 1472-2739 (on-line) 1472-2747 (printed)
\hfill {\pnum\folio}\else\ifodd\count0\def\\{ }%
\ifx\theshorttitle\relax\thetitle\else\theshorttitle\fi\hfill{\pnum\folio}
\else\def\\{ and }{\pnum\folio}\hfill\ifx\theshortauthors\relax\theauthors
\else\theshortauthors\fi\fi\fi}\vss}}
\footline{\vbox to 0pt{\vglue 0mm\line{\small\pfoot\ifnum\count0=\startpage
\copyright\ \gtp\hfill\else
\agt, Volume \thevolumenumber\ (\thevolumeyear)\hfill\fi}\vss}}
\else
\font=cmsl9 scaled 1050
\font\lnum=cmbx10 
\font=cmsl9 scaled 1050
\makeatletter
\def\@oddhead{{\small\lhead\ifnum\count0=\startpage ISSN 1472-2739 
(on-line) 1472-2747 (printed)\hfill {\lnum\number\count0}\else\ifodd\count0
\def\\{ }\ifx\theshorttitle\relax \thetitle \else\theshorttitle\fi\hfill
{\lnum\number\count0}\else\def\\{ and }{\lnum\number\count0}
\hfill\ifx\theshortauthors\relax 
\theauthors\else\theshortauthors\fi\fi\fi}}\def\@evenhead{\@oddhead}
\def\@oddfoot{\small\lfoot\ifnum\count0=\startpage\copyright\ \gtp\hfill\else
\agt, Volume \thevolumenumber\ (\thevolumeyear)\hfill\fi}
\def\@evenfoot{\@oddfoot}
\makeatother
\fi
\let\maketitlepage\makeagttitle
\let\makeshorttitle\maketitlepage
\let\maketitle\maketitlepage

\newwrite\gtoutfile
\long\gdef\makeheadfile{  
{\def\\{, }\def\s{ }
\immediate\openout\gtoutfile head.xxx
\immediate\write\gtoutfile{To: math@arxiv.org}
\immediate\write\gtoutfile{Subject: put OR rep NNNNN:ppppp}
\immediate\write\gtoutfile{--text follows this line--}
\immediate\write\gtoutfile{Proxy-for: \ifx\theasciiauthors\relax
\theauthors\else\theasciiauthors\fi\s<\ifx\theasciiemail\relax\theemail\else\theasciiemail\fi>}
\immediate\write\gtoutfile{\noexpand\\}
\immediate\write\gtoutfile{Authors: \ifx\theasciiauthors\relax
\theauthors\else\theasciiauthors\fi}
{\def\\{ }\immediate\write\gtoutfile{Title: \ifx\theasciititle\relax
\thetitle\else\theasciititle\fi}}
\immediate\write\gtoutfile{Subj-class: GT or SG, GR etc}
\immediate\write\gtoutfile{MSC-class: \theprimaryclass\ifx\thesecondaryclass\relax\else, \thesecondaryclass\fi}
\immediate\write\gtoutfile{Journal-ref: Algebr. Geom. Topol. \thevolumenumber\s
(\thevolumeyear) \startpage-\finishpage}
\immediate\write\gtoutfile{Comments: Published by Algebraic and
Geometric Topology at}
\immediate\write\gtoutfile{\s\s\s  http://www.maths.warwick.ac.uk/agt/AGTVol\thevolumenumber/agt-\thevolumenumber-\thepapernumber.abs.html}
\immediate\write\gtoutfile{\noexpand\\}
\immediate\write\gtoutfile{}
\ifx\theasciiabstract\relax
\immediate\write\gtoutfile{\theabstract}\else
\immediate\write\gtoutfile{\theasciiabstract}\fi
\immediate\write\gtoutfile{}
\immediate\write\gtoutfile{\noexpand\\}
\immediate\write\gtoutfile{}
\immediate\closeout\gtoutfile}}  

\def\maketitlepage{\makeagttitle\makeheadfile}
\let\makeshorttitle\maketitlepage
\let\maketitle\maketitlepage

\lognumber{10}
\volumenumber{2}
\volumeyear{2002}
\papernumber{10}
\published{27 March 2002}
\pagenumbers{219}{238}
\received{27 November 2001}
\accepted{25 February 2002}

\usepackage{amssymb,amsmath,amscd}

\newtheorem{thm}{Theorem}
\newtheorem{lem}{Lemma}
\newtheorem{prop}{Proposition}

\theoremstyle{definition}
\newtheorem{defn}{Definition} 
\newtheorem{exmp}{Example}

\newcommand{\Zset}{{\mathbb Z}}
\newcommand{\Cset}{{\mathbb C}}
\newcommand{\Rset}{{\mathbb R}}
\newcommand{\CP}{\Cset P}
\newcommand{\cpq}{\overline{\CP^2}}

\DeclareMathOperator{\Spc}{Spin^c}
\DeclareMathOperator{\Spin}{Spin}
\DeclareMathOperator{\SO}{SO}
\DeclareMathOperator{\im}{im}

\DeclareMathOperator{\ospc}{\Omega_4^{\Spc}}
\DeclareMathOperator{\tospc}{\tilde{\Omega}_4^{\Spc}}

\DeclareMathOperator{\tospin}{\tilde{\Omega}_4^{Spin}}
\DeclareMathOperator{\mspin}{MSpin}
\DeclareMathOperator{\mspc}{MSpin^{c}}
\DeclareMathOperator{\bspin}{BSpin}
\DeclareMathOperator{\bspc}{BSpin^{c}}
\DeclareMathOperator{\BSO}{BSO}
\DeclareMathOperator{\sign}{sign}
\DeclareMathOperator{\Arf}{Arf}
\DeclareMathOperator{\ks}{ks}
\begin{document}

\title{Stabilisation, bordism and embedded spheres\\in 4--manifolds}
\authors{Christian Bohr}                  
\address{Mathematisches Institut, Theresienstrasse 39\\80333 M\"unchen, 
Germany} 
\asciiaddress{Mathematisches Institut, Theresienstrasse 39\\80333 Muenchen, 
Germany} 
\email{bohr@mathematik.uni-muenchen.de}                     

\begin{abstract} 
It is one of the most important facts in 4--dimensional topology that
not every spherical homology class of a 4--manifold can be
represented by an embedded sphere. In 1978, M. Freedman and R. Kirby
showed that in the simply connected case, many of the obstructions to constructing such a sphere vanish if 
one modifies the ambient 4--manifold by adding products of \mbox{2--spheres}, a process which is
usually called stabilisation. In this paper, we extend this result to non--simply connected 4--manifolds
and show how it is related to the $\Spc$--bordism groups of Eilenberg--MacLane spaces.
\end{abstract}

\asciiabstract{ It is one of the most important facts in
4-dimensional topology that not every spherical homology class of a
4-manifold can be represented by an embedded sphere. In 1978,
M. Freedman and R. Kirby showed that in the simply connected case,
many of the obstructions to constructing such a sphere vanish if one
modifies the ambient 4-manifold by adding products of
2-spheres, a process which is usually called stabilisation. In
this paper, we extend this result to non-simply connected
4-manifolds and show how it is related to the Spin^c-bordism groups
of Eilenberg-MacLane spaces.}

\maketitle  

\section{Introduction and statement of results}

To determine the minimal genus of an embedded surface representing a given homology class
of a 4--manifold
has always been one of the most challenging problems in 4--dimensional topology.
In~\cite{KeM}, M. Kervaire and J. Milnor discovered the first non--trivial obstruction
to realising certain homology classes by embedded spheres, thus establishing that the
minimal genus is not always zero, not even in the simply connected case.
Their main result is the following.

\begin{thm}[Kervaire, Milnor]\label{kem}
Let $X$ be a closed, connected and oriented smooth 4--manifold and suppose that $\xi \in H_2(X;\Zset)$ is a characteristic class.
If $\xi$ can be represented by a smoothly embedded 2--sphere, then $\xi \cdot \xi \equiv \sign(X) \mod 16$.
\end{thm}

Recall that a homology class is called characteristic if the mod--2 reduction of its Poincar\'e dual
is the second Stiefel--Whitney class.
As an example, consider the case $X=\Cset P^2$ and let $\gamma$ denote the generator of the group 
$H_2(X;\Zset)=\Zset$ represented by a complex line.
It is not hard to see that the classes $\gamma$ and $2\gamma$ can both be represented by embedded spheres, 
whereas the obvious algebraic representative of $3\gamma$ 
has genus one. In fact,  Theorem \ref{kem} implies that 
this class cannot be represented by an embedded sphere.

Even if $X$ is simply connected, the converse of Theorem~\ref{kem} is of course not true. 
In the case $X=\CP^2$, P. Kron\-heimer and
T. Mrowka~\cite{KrM}  used Seiberg--Witten gauge theory to prove that the minimal genus in the class $n\gamma$ is 
given by $\frac{1}{2}(|n|-1)(|n|-2)$. 
Consequently $\gamma$ and $2\gamma$
are --- up to sign --- the
only classes  which can be realised by embedded spheres, although many other classes fulfill the
condition of Theorem~\ref{kem}.

However, it is known that all the gauge theoretical obstructions 
vanish if we pass from the manifold 
$X$ to one of the manifolds 
$X_k=X \# k(S^2 \times S^2)$ for large $k$. Note that
we have an inclusion $H_2(X;\Zset) \subset H_2(X_k;\Zset)$, 
so the statement that an embedded sphere in $X_k$
represents the class $\xi$ makes sense.
If we can find such a sphere for some $k$, we will say that the class $\xi$ can be
{\em stably represented} by an embedded sphere.
In 1978, Freedman and Kirby showed that if the 4--manifold $X$ is simply connected, 
the converse of Theorem \ref{kem} is true up to stabilisation.
More precisely, they proved the following result, which is Theorem~2 in~\cite{FK}.

\begin{thm}[Freedman, Kirby]\label{fk}
Let $X$ be a simply connected, closed and oriented smooth 4--manifold.
A characteristic homology class $\xi \in H_2(X;\Zset)$ 
can be stably represented by a smooth embedding of a 2--sphere  if and only if $\xi \cdot \xi \equiv \sign(X) \mod 16$.
\end{thm}

In this paper, we are interested in generalisations of this result to non--simply connected 4--manifolds.
First note that if $X$ is not simply connected, there may be homology classes which 
cannot even be represented by immersed spheres because they are not hit by the 
canonical map $\pi_2(X) \rightarrow H_2(X;\Zset)$. Homology classes in the image of this map are 
usually called spherical classes.
The main objective of this paper is to show that
Theorem~\ref{fk} is also true for 4--manifolds with non--trivial fundamental group,
provided that the classes in question are spherical, and to demonstrate that this fact
is a consequence of a non--trivial relation between $\Spin$ and $\Spc$--bordism groups
of Eilenberg--MacLane spaces.

\begin{thm}\label{FKlike}
Let $X$ be a closed, connected and oriented smooth 4--manifold.  
A characteristic spherical homology class $\xi \in H_2(X;\Zset)$ 
can be stably represented by a smoothly embedded sphere  if and only if $\xi \cdot \xi \equiv \sign(X) \mod 16$.
\end{thm}

Of course it is in general not true that a given immersion of a sphere is
homotopic to an embedding, not even stably. Instead there is a secondary obstruction which is defined
if a certain primary obstruction --- Wall's self intersection number --- vanishes and which determines 
whether a given homotopy class can be stably represented
by an embedding~\cite{ST}. The assertion of Theorem~\ref{FKlike} is that we can change the homotopy class
while fixing the homology class in such a way that all these obstructions vanish.

The technique applied in~\cite{FK} to construct an embedded sphere is to turn an
embedded surface of positive genus into a sphere by removing handles.
To be able to use this approach in the presence of a non--trivial fundamental
group, we have to find an embedding which enjoys the additional property
that the induced map from the fundamental group of the surface to the fundamental group of 
the ambient 4--manifold is trivial.

\begin{defn} An embedding $F \rightarrow X$ of a surface in a 
4--manifold $X$ will be called
{\em $\pi_1$--null} if the induced map $\pi_1(F) \rightarrow \pi_1(X)$ 
is trivial.
We will say that  a  homology class $\xi \in H_2(X;\Zset)$ can be 
{\em stably represented by a $\pi_1$--null embedding} 
if there is, for some natural number $k$, a $\pi_1$--null embedding $F \rightarrow X \# k(S^2 \times S^2)$
representing the class $\xi$.
\end{defn}

It is clear that a characteristic class which can be stably represented by a $\pi_1$--null embedding
is spherical. We shall see that this necessary condition is also sufficient, i.e.\ we have the following

\begin{thm}\label{charcase}
A characteristic homology class of a closed, connected and oriented smooth 4--manifold can be stably represented by
a smooth $\pi_1$--null embedding if and only if it is spherical.
\end{thm}

It should be mentioned that the condition that the homology class be characteristic is essential. In fact, by developing further
the arguments used in this paper, one can show~\cite{B} that there are spherical homology classes which are not characteristic and 
cannot be stably represented by a $\pi_1$--null embedding. In particular, they cannot
be stably represented by an embedded sphere, although they fulfill the condition 
$\xi \cdot \xi \equiv \sigma(X) \mod 16$. 
The reason is that there are  non--trivial 
additional obstructions to stably representing spherical classes by $\pi_1$--null embeddings. 
These obstructions are defined in terms of
certain bordism groups and are studied in~\cite{B}. The relation between stable embeddings of surfaces and
bordism theory which is investigated here and in~\cite{B}
is a generalisation of the ideas appearing in~\cite{KeM}. 
Note that the proof of Theorem~\ref{kem} is based on
Rokhlin's Theorem, which describes the image of the 4--dimensional $\Spin$--bordism group in the 4--dimensional $\Spc$--bordism group.

The proof of Theorem~\ref{charcase} will form the main part of this paper.
The idea of the proof is to use the language of B--structures and Kreck's stable diffeomorphism classification~\cite{K}
to reduce to the case that the class $\xi$ is trivial. This is done by splitting the manifold stably into a simply
connected 4--manifold and a spin manifold such that the homology class $\xi$ is moved into the simply connected part.
A crucial point is the choice of an appropriate bordism theory.
As the class $\xi$ is characteristic, it can be realised as the first Chern class of a $\Spc$--structure. This suggests
using $\Spc$--bordism groups so that we can keep track of the class $\xi$ by keeping track of the $\Spc$--structure.

In general, a splitting as indicated above can only be obtained after adding a simply connected 4--manifold,
therefore it is necessary to understand the effect of this on the representability by $\pi_1$--null embeddings.
To this end, we introduce a modified version of the
usual self intersection numbers of immersions in Section~\ref{constructing}, which may also
be of a certain interest in its own right.
The bordism problem into which the existence of the desired splitting can be translated 
is discussed in Section~\ref{bordism}. 
In Section~\ref{proofs}, we prove Theorem~\ref{charcase} and Theorem~\ref{FKlike},
and the last section of this manuscript is devoted to a short discussion of the topological versions
of our results.

Throughout this paper, all manifolds will be understood to be closed, connected, oriented and smooth,
unless stated otherwise.
An exception will be made in Section~\ref{top}, where we comment on the topological case.

The author would like to thank B. Hanke, D. Kotschick and R. Lee for helpful discussions
and reading preliminary versions of this manuscript and M. Kreck for explaining him some details
in~\cite{K}. He is also grateful to Yale University for its hospitality 
and to the {\it Deutsche Forschungsgemeinschaft} for financial support.

\section{Stable self intersection numbers}\label{constructing}

In this section, we introduce a modified version of Wall's  self intersection number which is 
adapted to our purposes. This modified invariant turns out to be the only obstruction to stably representing
spherical homology classes by $\pi_1$--null embeddings.
We will use this fact to show that in order for a spherical class to be stably representable by a $\pi_1$--null embedding,
it is sufficient to become representable after adding some simply connected 4--manifold, an observation which will
be used in the proof of Theorem~\ref{charcase}.

The relation between $\pi_1$--null embeddings and self intersection numbers of
immersions is provided by the following well known fact.

\begin{lem}\label{pi1null}
A homology class of a 4--manifold $X$
can be represented by a $\pi_1$--null embedding if and only if it
can be represented by an immersion $f \co S^2 \rightarrow X$ whose 
reduced self intersection number 
\[
\bar{\mu}(f) \in \Zset[\pi_1(X)] / (\langle \alpha - \alpha^{-1} \rangle + \Zset)
\]
vanishes.
\end{lem}

We refer the reader to~\cite{FQ} or~\cite{Wa} for a proof of this lemma as
well as for the definition and the most important
properties of the invariant $\bar{\mu}$.
For our purposes, it will be convenient to introduce a modified version of the
self intersection number which uses the homology classes of the loops associated with the self
intersection points instead of their homotopy classes.

\begin{defn}\label{defnofhomself}
Let $X$ be a 4--manifold and $f \co S^2 \rightarrow X$ an immersion. For a self intersection
point $x$ of $f$, we denote the associated group element by $c_x$.
The sum
\[
\mu_s(f) = \sum_x [ c_x ] \in H_1(X;\Zset_2)
\]
is called the {\em stable self intersection number} of $f$.
\end{defn}

Note that although $c_x \in \pi_1(X)$ is only defined up to orientation, its $\Zset_2$--homology class is well defined
and therefore the above expression makes sense.
The reason for the terminology will become clear later on as we shall see that one can modify the
reduced self intersection number of an immersion by stabilising, but that all these modifications do not
alter the stable self intersection number. 
For another description of the stable self intersection number
consider the canonical map $\Zset[\pi_1(X)] \rightarrow H_1(X;\Zset_2)$, given
by mapping $g \in \pi_1(X)$ to its $\Zset_2$--homology class. Note that an element of the
type $\alpha - \alpha^{-1}$ is mapped to $2\alpha=0$ and $1$ is mapped to $0$. Hence this map
induces a homomorphism
\[
R \co \Zset[\pi_1(X)] / (\langle \alpha - \alpha^{-1} \rangle + \Zset) \longrightarrow H_1(X;\Zset_2).
\]
It is immediate from the definition that, for an immersion $f$, the stable 
self intersection number $\mu_s(f)$ is simply the image of $\bar{\mu}(f)$ under the
above homomorphism. In particular, it only depends on the homotopy class of $f$
and defines a map $\mu_s\co \pi_2(X) \rightarrow H_1(X;\Zset_2)$,
which is the composition  $R \circ \bar{\mu}$.
Sometimes we will also use that $R$ factors over 
$\Zset_2[\pi_1(X)] / (\langle \alpha - \alpha^{-1} \rangle + \Zset_2)$.

\begin{exmp}\label{z2example}
Let us consider an easy example where we can actually compute the map $\mu_s$.
Consider the orientation preserving involution $t$ on $S^2 \times S^2$ given 
by $t(x,y)=(-x,\rho(y))$, where $\rho$ denotes the reflection at the hyperplane $x_3=0$.
Let $X = (S^2 \times S^2) / t$ and denote the projection by $\pi \co S^2 \times S^2 \rightarrow X$.
Clearly $\pi_1(X)=\Zset_2$. The two spheres $\alpha = S^2 \times 1$ and
$\beta=1 \times S^2$ define a hyperbolic basis $\alpha, \beta \in H_2(S^2 \times S^2;\Zset)$.
A straightforward calculation shows that $\pi_*(\beta)=0$
whereas $\pi_*(\alpha)$ is the non--zero element of $H_2(X;\Zset)=\Zset_2$, and that the action
of the deck transformation group on $H_2(S^2 \times S^2;\Zset)$ is just the multiplication by $-1$. 
Hence we obtain
\begin{align*}
\mu_s(x\alpha + y \beta) & =  xy  \mod 2 \\
\pi_*(x\alpha + y\beta) &= x  \mod 2.
\end{align*}
Since the classes $\alpha$ and $\alpha+\beta$ both
project to the non--trivial element of $H_2(X;\Zset)$ but have different stable self intersection 
numbers,
this example shows in particular that the stable self intersection number of an immersion 
depends in general really on its homotopy class and not only on its homology class.
\end{exmp}

The following proposition is the main technical result of this section.
It shows that the relation between stable self intersection numbers and stable $\pi_1$--null
embeddings is similar to that between Wall's self intersection number and $\pi_1$--null embeddings.

\begin{prop}\label{embeddings}
Let $X$ be a 4--manifold and $\xi \in H_2(X;\Zset)$ a 
homology class. Then the following conditions
are equivalent.
\begin{enumerate}
\item\label{condi} 
There is an immersion $f \co S^2 \rightarrow X$ representing $\xi$ such 
that $\mu_s(f)=0$.
\item\label{condii} 
There is, for some $k$, a $\pi_1$--null embedding $F \rightarrow X \# k(S^2 \times S^2)$
representing the homology class $(\xi,0, \dots, 0)$.
\item\label{condiii} 
There is, for some simply connected 4--manifold $Y$, a $\pi_1$--null embedding $F \rightarrow
X \# Y$ such that the homology class of $F$ is $(\xi,c)$ with a characteristic class $c \in H_2(Y;\Zset)$.
\end{enumerate}
\end{prop}

Before we can start with the proof of this result, we have to state and prove some technical
lemmas, the first one being a simple algebraic observation. 

\begin{lem}\label{observation}
Suppose that $Q \co \Gamma \times \Gamma \rightarrow \Zset$ is a symmetric unimodular bilinear form
over the integers, where $\Gamma$ is a free $\Zset$--module,
and $x,y \in \Gamma$ are elements whose sum $x+y$ is characteristic. Then
$Q(x,y) \equiv 0 \mod 2$.
\end{lem}

\begin{proof}
Since $x+y$ is characteristic, we have
\[
Q(x,x) \equiv Q(x+y,x) \equiv Q(x,x) + Q(x,y) \mod 2,
\]
and the assertion follows.  
\end{proof}

\begin{lem}\label{stablered}
Suppose that $X$ is a 4--manifold, 
$[f] \in \pi_2(X)$, $g \in \pi_1(X)$, and $\epsilon \in \{-1,+1\}$.
Then there exists a homotopy class $[f'] \in \pi_2(X \# (S^2 \times S^2))$ such that
\begin{enumerate}
\item $f'_*[S^2]=(f_*[S^2],0) \in H_2(X \# (S^2 \times S^2);\Zset)$, and
\item $\bar{\mu}(f')=\bar{\mu}(f)+2\epsilon g$.
\end{enumerate}
\end{lem}

\begin{proof} 
Let $Y=X \# (S^2 \times S^2)$ and denote by $\tilde{Y}$ the universal covering. 
This manifold is obtained from the universal covering $\tilde{X}$ of $X$ by adding 
copies of $S^2 \times S^2$, one at each preimage of the point at which we performed
the connected sum of $X$ and $S^2 \times S^2$.
It follows easily from the
obvious Mayer--Vietoris sequence that we have a splitting
\[
\pi_2(Y)=H_2(\tilde{Y};\Zset) = H_2(\tilde{X};\Zset) \oplus (H_2(S^2 \times S^2;\Zset) 
\otimes_{\Zset} \Zset[\pi_1(X)])
\]
which is orthogonal with respect to the intersection form. The action of the deck transformation group
on the second summand (from the left) is given by $g(x \otimes h)=(x \otimes gh)$, and 
the restriction of the intersection form $\lambda$ with values in the group ring to
the second summand is completely described  by $\lambda((0,x \otimes g),(0,y \otimes h))=(x \cdot y)gh^{-1}$,
where the dot denotes the intersection product on $H_2(S^2 \times S^2;\Zset)$.

Now choose some $x \in H_2(S^2 \times S^2;\Zset)$ such that $x \cdot x = -2\epsilon$. Consider the element
$\beta = (0,x \otimes 1 - x \otimes g) \in H_2(\tilde{Y};\Zset)$.
An easy computation shows that $\lambda(\beta,\beta)=-4 \epsilon + 2 \epsilon (g + g^{-1})$ and 
$\bar{\mu}(\beta)=  2 \epsilon g$. If $\pi \co H_2(\tilde{Y};\Zset) \rightarrow H_2(Y;\Zset)$ 
denotes the projection, we
have $\pi_*(\beta)=0$.

The element $[f] \in \pi_2(X) \subset \pi_2(Y)$ is contained in the first summand of the above decomposition.
Choose a map $f' \co S^2 \rightarrow Y$ in the homotopy class $[f] + \beta$.
Then $f'_*[S^2]=\pi_*([f]+\beta)=\pi_*([f])=(f_*[S^2],0)$, and, since $\lambda([f],\beta)=0$, we have
$\bar{\mu}(f')=\bar{\mu}(f) + 2 \epsilon g$. 
\end{proof}

\begin{lem}\label{boundaries}
Suppose that $X$ is a 4--manifold and $[f] \in \pi_2(X)$, $g,h \in \pi_1(X)$. Then there is a homotopy class
$[f'] \in \pi_2(X \# (S^2 \times S^2))$ such that
\begin{enumerate}
\item $f'_*[S^2]=(f_*[S^2],0) \in H_2(X \# (S^2 \times S^2);\Zset)$, and
\item $\bar{\mu}(f')=\bar{\mu}(f)+g+h-gh$.
\end{enumerate}
\end{lem}

\begin{proof} 
Let $Y=X \# (S^2 \times S^2)$ and denote by $\tilde{Y}$ the universal covering. Again we will make use
of the decomposition
\[
\pi_2(Y) = H_2(\tilde{Y};\Zset) = H_2(\tilde{X};\Zset) \oplus 
(H_2(S^2 \times S^2;\Zset) \otimes_{\Zset} \Zset[\pi_1(X)]).
\]
As in the proof of Lemma~\ref{stablered}, we only have to find an element $\beta$ in the
second summand of this decomposition such that $\pi_*(\beta)=0$ and $\bar{\mu}(\beta)=g+h-gh$.
For this purpose, consider the elements
\begin{align*}
x &= (1,-1) \\
y &= (0,1) \\
z &= (-1,0)
\end{align*}
in the group $H_2(S^2 \times S^2;\Zset)$, here we use, as usual, the basis of this group given by the
homology classes of the factors. Then $x \cdot y = x \cdot z =1$, $ y \cdot z = -1$ and
$x + y + z=0$. Now let
\[
\beta= (0, x\otimes 1 + y \otimes g + z \otimes h^{-1}).
\]
Then $\pi_*(\beta)=x+y+z=0$. Furthermore, a short calculation shows that $\bar{\mu}(\beta)=g+h-gh$,
and the proof is complete.  
\end{proof}

\begin{proof}[Proof of Proposition~\ref{embeddings}]
(\ref{condi}) $\implies$ (\ref{condii}): Assume that we are given an immersion $f$ as in the statement of
the proposition.
Pick a representative of $\bar{\mu}(f)$ in the group ring $\Zset[\pi_1(X)]$ and let
$x \in \Zset_2[\pi_1(X)]$ denote its mod--2 reduction.
By assumption, the image of $x$ under the canonical map $\Zset_2[\pi_1(X)] \rightarrow H_1(X;\Zset_2)$
is zero. 

In  the non--homogeneous description of the standard chain complex 
of $\pi_1(X)$ with $\Zset_2$--coefficients, the
group $C_1(\pi_1(X);\Zset_2)$ of one--cycles consists of all finite linear combinations
$\sum_{g \in \pi_1(X)} n_g [g]$, $n_g \in \Zset_2$ (see for instance \cite{Br}, Chapter~II.3).
The boundary of an element $[g|h] \in C_2(\pi_1(X);\Zset_2)$ 
is given by 
\[
\partial[g|h]=[g]+[h]+[gh].
\] 
If we use the 
identification of $\Zset_2[\pi_1(X)]$ and $C_1(\pi_1(X);\Zset_2)$ 
as $\Zset_2$--modules, induced by $g \mapsto [g]$, the kernel of the
map $\Zset_2[\pi_1(X)] \rightarrow H_1(X;\Zset_2)$ is identified with the kernel of the projection 
\[
C_1(\pi_1(X);\Zset_2) \longrightarrow H_1(\pi_1(X);\Zset_2),
\] 
which is by definition the 
subgroup generated by all boundaries $[g]+[h]+[gh]$. 
Consequently, the element $x$ is a sum $x=\sum_j (g_j+h_j+g_j h_j)$ for some elements 
$g_j,h_j \in \pi_1(X)$. 
By Lemma~\ref{boundaries}, we can therefore construct an element 
$[f'] \in  \pi_2(X \# k(S^2 \times S^2))$ such that
\[
\bar{\mu}(f') \equiv \bar{\mu}(f) + \sum_j (g_j + h_j - g_j h_j) \equiv x + x \equiv 0 \mod 2,
\]
and $f'_*[S^2]=(f_*[S^2],0,\dots,0)$. 
This implies that there is a representative $\sum_i a_i g_i$ of
$\bar{\mu}(f')$ in $\Zset[\pi_1(X)]$ such that all the coefficients $a_i$ are even numbers. By 
Lemma~\ref{stablered}, we can now find an immersion $f'' \co S^2 \rightarrow X \# k'(S^2 \times S^2)$
for some $k'$ which represents $(f_*[S^2],0, \dots, 0)$ and has reduced self intersection number zero. 
The result then follows from Lemma~\ref{pi1null}.

\noindent (\ref{condii}) $\implies$ (\ref{condiii}): clear, since $S^2 \times S^2$ is spin, so the class 
$0 \in H_2(S^2 \times S^2;\Zset)$
is characteristic.

\noindent (\ref{condiii}) $\implies$ (\ref{condi}): Assume that we are given a $\pi_1$--null 
embedding $\iota \co F \rightarrow X\# Y$ representing a class $(\xi,c)$, where $c$ is characteristic. 
Let $Z=X \# Y$ and denote by $\pi \co \tilde{Z} \rightarrow Z$ the universal covering.
By Lemma~\ref{pi1null}, there is an immersion $f \co S^2 \rightarrow Z$ representing $(\xi,c)$
such that $\bar{\mu}(f)=0$. Let $[f] \in H_2(\tilde{Z};\Zset)=\pi_2(Z)$ denote the corresponding
homology class. Then $\pi_*([f])=(\xi,c) \in H_2(Z;\Zset)$.
We have a decomposition
\[
H_2(\tilde{Z};\Zset) = H_2(\tilde{X};\Zset) \oplus (H_2(Y;\Zset) \otimes_{\Zset} \Zset[\pi_1(X)]).
\]
Let $[f]=\alpha +\beta$ be the corresponding decomposition of the element $[f]$. 
Using that $\lambda(\alpha,\beta)=0$, it is easy to see 
that 
$0=\mu_s([f])=\mu_s(\alpha)+\mu_s(\beta)$, and clearly $\pi_*(\alpha)=\xi$, $\pi_*(\beta)=c$.
If we can prove that $\mu_s(\beta)=0$,
we obtain $\mu_s(\alpha)=0$ and any immersion corresponding to $\alpha$ will fulfill our requirements.

In order to show that $\mu_s(\beta)=0$, write 
$\beta=\sum_{i=1}^n x_i \otimes g_i \in H_2(Y;\Zset) \otimes_{\Zset} \Zset[\pi_1(X)]$,
where $x_i \in H_2(Y;\Zset)$ and $g_i \in \pi_1(X)$.
An easy calculation shows that
$\bar{\mu}(\beta)= \sum_{i < j} (x_i \cdot x_j) g_i g_j^{-1}$.
Hence we have
\begin{align*}
\mu_s(\beta) &= \sum_{i < j} (x_i \cdot x_j) ([g_i] + [g_j]) \\
&= \sum_{i < j} (x_i \cdot x_j) [g_i] + \sum_{j < i} (x_i \cdot x_j) [g_i] \\
&= \sum_i \Big ( \sum_{j > i} x_i \cdot x_j \Big ) [g_i] + \sum_i \Big ( \sum_{j < i} x_i \cdot x_j \Big ) [g_i] \\
&= \sum_i \Big ( \sum _{j \neq i} x_i \cdot x_j \Big ) [g_i].
\end{align*}
Now we know that $\pi_*(\beta)=c$ is characteristic, i.e $\sum_i x_i$ is a characteristic class for $Y$. Using Lemma 
\ref{observation}, we therefore obtain (for fixed $i$)
\[
\sum_{j \neq i} x_i \cdot x_j = x_i \cdot (\sum_{j \neq i} x_j) = x_i \cdot (c-x_i) \equiv 0 \mod 2,
\]
since $x_i + (c-x_i)=c$ is characteristic. As this is true for every $i$, 
all the coefficients in the above expression are 
even and therefore $\mu_s(\beta)=0$ as desired.   
\end{proof}

\section{Some remarks on Spin${}^c$--bordism groups}\label{bordism}

In this section, we shall prove a result
about $\Spc$--bordism groups (Proposition~\ref{exactness}) of which Theorem~\ref{charcase} will be a consequence.
First let us summarize some well known facts about
$\Spin$ and $\Spc$--bordism groups, which can for instance be found in \cite{St}.

\begin{lem}
\mbox{}
\begin{enumerate}
\item  $\Omega_0^{\Spc}(*)=\Zset$, $\Omega_1^{\Spc}(*)=\Omega_3^{\Spc}(*)=0$ and 
$\Omega_2^{\Spc}(*)=\Zset$.
\item $\Omega_0^{\Spin}(*)=\Zset$, $\Omega_1^{\Spin}(*)=\Omega_2^{\Spin}(*)=\Zset_2$ and 
$\Omega_3^{\Spin}(*)=0$.
\end{enumerate}
\end{lem}

Throughout this section, all CW--complexes will be assumed to be connected and to have only one 0--cell,
which we use as base point.

\begin{defn} Let $K$ be a space.
Define a homomorphism
\[
c_1 \co \ospc(K) \longrightarrow H_2(K;\Zset)
\]
by $c_1([(X,s,f)])=f_*PD(c_1(s))$, where $X$ is a connected 4--mani\-fold with a normal $\Spc$--structure
$s$, $c_1(s)$ is the first Chern class of $s$ and $f \co X \rightarrow K$ is a 
continuous map.
\end{defn}

This map is clearly well--defined and natural in $K$.
For the remainder of this section, we will work with the reduced homology theory $\tilde{\Omega}_*^{\Spc}$
instead of $\Omega_*^{\Spc}$.
Since $c_1$
clearly vanishes on the bordism group of the point, it defines a homomorphism
\[
c_1 \co \tilde{\Omega}_4^{\Spc}(K) \longrightarrow H_2(K;\Zset).
\]
Our next aim will be to describe this homomorphism in greater detail and to understand its kernel.
For this purpose, let us consider the 
Atiyah--Hirzebruch spectral sequence (AHSS for short) 
\[
E^2_{p,q}=\tilde{H}_p(K;\tilde{\Omega}_q^{\Spc}(S^0)) \Longrightarrow \tilde{\Omega}_{p+q}^{\Spc}(K).
\]
Clearly $E^{\infty}_{0,4}=E_{3,1}^{\infty}=E_{1,3}^{\infty}=0$.
As $\Omega_4^{\SO}(K) \rightarrow H_4(K;\Zset)$ is onto and every orientable 4--manifold has
a normal $\Spc$--structure,
all differentials emerging from $E^*_{4,0}$ vanish.  
Hence we have a short exact sequence
\[
0 \longrightarrow E^{\infty}_{2,2} \longrightarrow \tospc(K) \longrightarrow H_4(K;\Zset) \longrightarrow 0.
\] 
Since all differentials emerging from $E_{2,2}^*$ are zero,
we have a natural projection $E_{2,2}^2=H_2(K;\Zset) \rightarrow E_{2,2}^{\infty}$, 
and we obtain a natural map 
\[
H_2(K;\Zset) \longrightarrow \tospc(K),
\] 
whose kernel is the image of the differential
$d_3 \co H_5(K;\Zset) \rightarrow H_2(K;\Zset)$.

\begin{defn} For every CW--complex $K$ let $\varphi$ denote the composition
\[
\varphi \co H_2(K;\Zset) \longrightarrow \tospc(K) \stackrel{c_1}{\longrightarrow} H_2(K;\Zset)
\]
of the map described above and $c_1$.
\end{defn}

\begin{lem}\label{varphi} 
The homomorphism $\varphi$ is the multiplication by $2$.
\end{lem}

\begin{proof}
First, we show this in the special case that $K$ is a connected orientable surface.
Then $H_2(K;\Zset)=\Zset$, so $\varphi \co H_2(K;\Zset) \rightarrow H_2(K;\Zset)$ is the multiplication
by some number $n$. We have to prove that $n=2$.

For this purpose, suppose that $[(X,s,f)]$ is some element in $\tospc(K)$. 
Since $c_1(s) \equiv w_2(X) \mod 2$, 
we have $c_1([(X,s,f)]) \equiv f_*PD(w_2(X)) \mod 2$. 
Now it is an immediate consequence of Wu's Theorem that the second
Stiefel--Whitney class of a 4--manifold and its orientation class are related by
\[
PD(w_2(X))=(Sq^2)^* [X]_2,
\] 
where $(Sq^2)^*$ denotes the homology operation dual to the second
Steenrod square and $[X]_2$ is the mod--2 orientation.  By naturality, we obtain that
\[
c_1([(X,s,f)]) \equiv f_*((Sq^2)^* [X]_2)=(Sq^2)^*f_*[X]_2=0 \mod 2,
\]
since $H_4(K;\Zset_2)=0$. This implies that every element in the image of $\varphi$ is divisible by $2$, 
hence $n$ must be even.

So we are done if we can show that every multiple of $2$ is in the image of $\varphi$.
Suppose that $x=2y \in H_2(K;\Zset)$. Let $X=K \times K$ and let $f\co X \rightarrow K$ denote 
the projection to
one factor. Then, by the K\"unneth Theorem, $f_* \co H_2(X;\Zset) \rightarrow H_2(K;\Zset)$ is onto. 
Pick a preimage
$\bar{y}$ of $y$ and let $\bar{x}=2\bar{y}$. Since $X$ is spin
(and hence the stable normal bundle is spin), there exists a normal $\Spc$--structure $s$ on
$X$ having first Chern class $c_1(s)=PD(\bar{x})$, hence $c_1([(X,s,f)])=x$. Moreover $H_4(K;\Zset)$=0, so
$E_{2,2}^{\infty}=E_{2,2}^2=H_2(K;\Zset)$, 
and $[(X,s,f)]$ defines an element $\alpha \in E^{\infty}_{2,2}=H_2(K;\Zset)$.  By definition,
$\varphi(\alpha)=c_1([(X,s,f)])=x$.
So we have proved that $x$ is in the image of $\varphi$, and this
completes the proof of the lemma in the special case that $K$ is a surface.

As to the general case, note that for every $x \in H_2(K;\Zset)$, there is a connected oriented surface $F$
and a map $f\co F \rightarrow K$ mapping the orientation of $F$ to $x$. 
Since $\varphi$ is the multiplication by
$2$ on $H_2(F;\Zset)$, we have
\[
\varphi(x)=\varphi(f_*[F])=f_*\varphi([F])=f_*(2[F])=2x,
\]
and this proves the lemma in the general case.  
\end{proof}

As for $\Spc$--bordism, we also have reduced $\Spin$--bordism groups,
and for each $CW$--complex $K$, we have a natural map
\[
\tospin(K) \longrightarrow \tospc(K)
\]
whose image is clearly in the kernel of $c_1$. To understand this image, we will need the following
facts about the AHSS
\[
\tilde{H}_p(K;\tilde{\Omega}_q^{\Spin}(S^0)) \Longrightarrow \tilde{\Omega}_{p+q}^{\Spin}(K).
\]
which are easily verified by
comparing with the AHSS for 
$\tilde{\Omega}_*^{\Spin}(\Rset P^{\infty})$
and using the computations of the differentials in the $E^2$--term carried out in~\cite{Tei}.

\begin{lem}\label{vanishing}
For every CW--complex $K$, 
the differential 
\[
d_3 \co E^3_{4,0} \longrightarrow E^3_{1,2}
\]
in the third term of the
Atiyah--Hirzebruch spectral sequence for 
$\tilde{\Omega}_*^{\Spin}(K)$ 
vanishes.
\end{lem}

\begin{lem}\label{image}
Suppose that $K$ is a CW--complex. The term $E_{4,0}^\infty$ in the AHSS for
$\tilde{\Omega}_*^{\Spin}(K)$ is exactly the kernel of 
\[
(Sq^2)^* \circ red_2 \co H_4(K;\Zset) \longrightarrow H_2(K;\Zset_2),
\]
where $(Sq^2)^*$ is the operation dual to $Sq^2$ and $red_2$ denotes reduction mod 2.
\end{lem}

To analyze the relation between $\tilde{\Omega}_*^{\Spin}$ and $\tilde{\Omega}_*^{\Spc}$, it proves
useful to study the ``relative bordism theory'' associated to the map $\bspin \rightarrow \bspc$.
Geometrically, the objects in this bordism theory can be described as 
manifolds with normal $\Spc$--structures and liftings
of these structures to $\Spin$--structures on the boundary.
However, we will use the description in terms of the associated spectra.
Let $\mspin$ and $\mspc$ denote the Thom spectra of
$\bspin$ respectively $\bspc$. There is an obvious map $\mspin \rightarrow \mspc$
which we can assume to be an inclusion. We denote the cone over $\mspin$ by $C(\mspin)$.

\begin{defn}
Let $G=\mspc \cup C(\mspin)$.
We will also denote the associated reduced homology theory by $G_*$, i.e.\ for
a space $K$ with base point, we have
\[
G_*(K)= \pi_*(G \wedge K).
\]
\end{defn}

The mappings $\mspin \rightarrow \mspc$ and $\mspc \rightarrow G$ of spectra induce natural transformations
$\tilde{\Omega}^{\Spin}_* \rightarrow \tilde{\Omega}_*^{\Spc}$ and 
$\tilde{\Omega}_*^{\Spc} \rightarrow G_*$ between the corresponding homology theories.
The following lemma summarizes some basic facts about
the homology theory $G_*$ and these natural transformations which are easily proved using the fact that
$\mspin \rightarrow \mspc \rightarrow G$ is a cofibre sequence and some standard results for which the reader
is referred to \cite{Ad} or \cite{Rd}.

\begin{lem}\label{G}
\mbox{}
\begin{enumerate}
\item For every CW--complex $K$, the sequence
\[
\tilde{\Omega}_*^{\Spin}(K) \longrightarrow \tilde{\Omega}_*^{\Spc}(K) \longrightarrow G_*(K) 
\]
is exact.
\item $G_2(S^0)=\Zset$, $G_3(S^0)=\Zset_2$ and $G_i(S^0)=0$ for $i \leq 1$.
\item The map $\tilde{\Omega}_2^{\Spc}(S^0) \rightarrow G_2(S^0)$ is the multiplication by $2$.
\end{enumerate}
\end{lem}

Using these results and the fact that the natural transformation $G_* \rightarrow \tilde{\Omega}_{*-1}^{\Spin}$ 
given by the map $G \rightarrow \Sigma\mspin$ induces maps between the corresponding spectral sequences,
we can now determine some differentials in the AHSS for $G_*$. Again we omit the proof which is straightforward.

\begin{lem}\label{d2zero}
Let $K$ be a CW--complex and consider the AHSS 
\[
E_{p,q}^2=\tilde{H}_p(K;G_q(S^0)) \Longrightarrow G_{p+q}(K)
\]
for $G_*(K)$. We then have 
\begin{enumerate}
\item 
For every $p \geq 3$, the differential 
\[
d_2 \co H_p(K;G_2(S^0)) \longrightarrow H_{p-2}(K;G_3(S^0))
\] 
is the reduction mod 2 followed by the operation $(Sq^2)^*$ dual to $Sq^2$.
\item 
$E^{\infty}_{1,3}=E^2_{1,3}=H_1(K;\Zset_2)$.
\item $E^{\infty}_{2,2}=E^2_{2,2}=H_2(K;\Zset)$.
\end{enumerate}
\end{lem}

\begin{lem}\label{intersection}
Suppose that $K$ is a CW--complex with 1--skeleton $K^1$. Then
\[
\im (G_4(K^1) \longrightarrow G_4(K)) \cap \im (\tospc(K) \longrightarrow G_4(K)) = \{ 0 \}.
\]
\end{lem}

\begin{proof}
Suppose we are given some $\alpha \in \im(G_4(K^1) \rightarrow G_4(K))$ and some 
$\beta \in \tospc(K)$ mapping to $\alpha$ under $\tospc(K) \rightarrow G_4(K)$. By definition,
the subgroup $\im G_4(K^1)$ is just the first stage $F_1G_4(K)$ of the filtration of $G_4(K)$ used in the 
AHSS converging to $G_*(K)$, and $F_1G_4(K)=E^{\infty}_{1,3}$.

Now suppose that $\alpha \neq 0 \in E^{\infty}_{1,3}$. Recall that 
$E^{\infty}_{1,3}=E^2_{1,3}=H_1(K;\Zset_2)$
by Lemma~\ref{d2zero}. We can find a map $f \co K \rightarrow \Rset P^{\infty}$ such that
$f_*\alpha \neq 0$, and $f$ can be chosen to be cellular. 
Let $\bar{E}^*_{*,*}$ denote the terms of the AHSS for $G_*(\Rset P^{\infty})$. Again we have
$\bar{E}^{\infty}_{1,3}=\bar{E}^2_{1,3}=H_1(\Rset P^{\infty};\Zset_2)$, so 
$f_* \co E^{\infty}_{1,3} \rightarrow \bar{E}_{1,3}^{\infty}$ maps $\alpha$ to the non--zero element
of $\bar{E}_{1,3}^{\infty}=\Zset_2$.

Now $\bar{E}^{\infty}_{1,3}=F_1G_4(\Rset P^{\infty})$, 
hence we obtain that 
$f_*\alpha \neq 0 \in G_4(\Rset P^{\infty})$. On the other hand,
we have a commuting diagram
\[
\begin{CD}
\tospc(K) @>{f_*}>> \tospc(\Rset P^{\infty}) \\
@VVV @VVV \\
G_4(K) @>{f_*}>> G_4(\Rset P^{\infty}) 
\end{CD}
\]
and as $H_2(\Rset P^{\infty};\Zset)=H_4(\Rset P^{\infty};\Zset)=0$, we have
$\tospc(\Rset P^{\infty})=0$, in particular $f_*\beta=0$. This shows that $f_*\alpha=0$, 
a contradiction. 
\end{proof}

\begin{lem}\label{specialcase} 
Suppose that $K$ is a CW--complex and 
that $\alpha \in \tospc(K)$. If the image of $\alpha$ in
$H_4(K,\Zset)$ is zero and
$c_1(\alpha)=0$, then $\alpha$ is in the
image of the natural map $\tospin(K) \rightarrow \tospc(K)$.
\end{lem}

\begin{proof}
By Lemma \ref{G}, we only have to show that the image
of $\alpha$ under the map $\tospc(K) \rightarrow G_4(K)$ is zero.

Let $E_{*,*}^*$ denote the terms of the AHSS for $\tilde{\Omega}_*^{\Spc}(K)$ and denote by 
$\bar{E}_{*,*}^*$ 
the terms of the spectral sequence for $G_*(K)$. The CW--decompo\-sition of $K$ yields filtrations 
$F_*\tospc(K)$ of $\tospc(K)$ and $F_*G_4(K)$ of $G_4(K)$, and we have a
filtration preserving homomorphism $\tospc(K) \rightarrow G_4(K)$. 
As the image of $\alpha$ in $H_4(K;\Zset)$ is zero, $\alpha \in F_2\ospc(K)$ 
(recall that $E^{\infty}_{3,1}=0$).
Let $\bar{\alpha} \in F_2G_4(K)$ denote the image of $\alpha$. We claim that $\bar{\alpha} \in F_1G_4(K)$,
i.e.\  that the induced map $E_{2,2}^{\infty} \rightarrow \bar{E}_{2,2}^{\infty}$ maps the equivalence class
$[\alpha]$ of $\alpha$ to zero.
In fact, $E^{\infty}_{2,2}$ is a quotient of $E^2_{2,2}$ and, by Lemma~\ref{d2zero}, 
$\bar{E}^{\infty}_{2,2}=E^2_{2,2}$, so we have a commuting diagram:
\[
\begin{CD}
E^2_{2,2} @>>> \bar{E}_{2,2}^2 \\
@VVV @VVV \\
E_{2,2}^{\infty} @>>> \bar{E}_{2,2}^{\infty}
\end{CD}
\]
Pick a lift $x \in E^2_{2,2}=H_2(K;\Zset)$ of $[\alpha]$. According to the definition of $\varphi$ and by
Lemma~\ref{varphi},
we have $\varphi(x)=c_1(\alpha)=0=2x$, so $x$ has order at most $2$.
Now Lemma~\ref{G} implies that the upper horizontal arrow is the multiplication by~$2$.
Hence it maps $x$ to zero and we obtain that
$\bar{\alpha} \in F_1G_4(K)$, as claimed.
 
Now, by definition, $F_1G_4(K)$ is simply the image of $G_4(K^1) \rightarrow G_4(K)$, where $K$ denotes
the 1--skeleton of $K$. Hence, by Lemma \ref{intersection}, we can conclude that actually
$\bar{\alpha}=0$, and the lemma is proved.   
\end{proof}

\begin{prop}\label{exactness}
For every CW--complex $K$, the sequence
\[
\tospin(K) \longrightarrow \tospc(K) \stackrel{c_1}{\longrightarrow} H_2(K;\Zset)
\]
is exact.
\end{prop}

\begin{proof}
It is immediate from the definition of $c_1$ that we have the inclusion image $\subset$ kernel.
So assume that we have a class $\alpha \in \tospc(K)$ such that $c_1(\alpha)=0$. Let $(X,s,f)$ be
a triple representing $\alpha$. As in the proof of Lemma \ref{varphi}, we then have
\[
(Sq^2)^*f_*[X]_2 = f_*PD(w_2)= red_2(c_1(\alpha))=0 \in H_2(K;\Zset_2).
\]
But by Lemma \ref{image},
the image of $\tospin(K)$ in $H_4(K;\Zset)$ is precisely the kernel of $(Sq^2)^* \circ red_2$, and therefore
we can find an element $\beta \in \tospin(K)$ having the same image in $H_4(K;\Zset)$ as $\alpha$. Let
$\bar{\beta}$ denote the image of $\beta$ in $\tospc(K)$. Then $\alpha-\bar{\beta}$ is in the kernel of
$\tospc(K) \rightarrow H_4(K;\Zset)$ and $c_1(\alpha-\bar{\beta})=c_1(\alpha)=0$. 
By Lemma~\ref{specialcase},
we can conclude that $\alpha-\bar{\beta}$ 
is in the image of $\tospin(K)$, and since $\bar{\beta} \in \im \tospin(K)$
by construction, we are done.  
\end{proof}

\section{The proofs of Theorem~\ref{charcase} and Theorem~\ref{FKlike}}\label{proofs}

In this section, we will use Proposition~\ref{exactness} to prove Theorem~\ref{charcase}
and Theorem~\ref{FKlike}. We mention that, at least in the case that $K$ is an Eilenberg--MacLane space,
Proposition~\ref{exactness} is in fact equivalent to the statement of Theorem~\ref{charcase}. The reason is that,
after adding copies of $\CP^2$, one can turn a stable $\pi_1$--null embedding realising a characteristic class into
an embedding of a sphere along which one can attach a 3--handle to obtain a spin manifold.

\begin{proof}[Proof of Theorem~\ref{charcase}]
Since a $\pi_1$--null embedding can be lifted to the universal covering, it is clear that
every class which can be stably represented by a $\pi_1$--null embedding is spherical.
To prove the converse,
assume that we are given a 4--manifold $X$ and a spherical and characteristic class $\xi \in H_2(X;\Zset)$.

First let us assume that there is a spherical class $\omega$ such that $\omega \cdot \xi = 1$.
We will use the abbreviation $\Pi=\pi_1(X)$. Pick a map 
$u \co X \rightarrow K(\Pi,1)$ inducing an isomorphism
$\pi_1(X) \rightarrow \Pi$ and choose a normal $\Spc$--structure $s$ on $X$ having 
first Chern class $c_1(s)=PD(\xi)$.
By definition, the homomorphism $c_1$ maps the bordism class $(X,s,u)$ to $u_*\xi$. As $\xi$ is spherical,
this is zero. Hence Proposition~\ref{exactness} yields a $\Spc$ 4--manifold $(Z,s_1)$ and a spin 4--manifold $Y$
together with a map $v \co Y \rightarrow K(\Pi,1)$ such that
\[
[(Z,s_1,*)] + [(Y,s_2,v)] = [(X,s,u)] \in \Omega_4^{\Spc}(K(\Pi,1)).
\]
Here $s_2$ denotes the $\Spc$--structure induced by the spin structure on $Y$. Note
that $c_1(s_2)=0$. We can also assume that $Z$
is simply connected and that $v_* \co \pi_1(Y) \rightarrow \Pi$ is an isomorphism. Furthermore we can arrange for
$c_1(s_1)$ to be a primitive class by adding a copy of $\CP^2 \#\cpq$ with the unique $\Spc$--structure
having Chern class $(1,1)$.

Let $B=\bspc \times K(\Pi,1)$ and consider the obvious fibration $B \rightarrow \BSO$. 
The normal $\Spc$--structure $s$ and the map $u$ define a $B$--structure $\nu \co X \rightarrow B$ on $X$.
It is not difficult to see that the existence of the class $\omega$ implies that $\nu$ is a 2--equivalence.
Similarly the $\Spc$--structure $s'$ on $X'=Z \# Y$ obtained by gluing $s_1$ and $s_2$
and the map $u' \co X' \rightarrow \Pi$ which is $v$ on $Y$ and trivial on $Z$ define a $B$--structure 
$\nu' \co X' \rightarrow B$. As the class $c_1(s')=c_1(s_1)$ is primitive and 
$\nu'_* \co \pi_1(X') \rightarrow \Pi$ is an isomorphism, this $B$--structure is also a 2--equivalence.

Of course the bordism group of $B$--structures is simply $\Omega_*^{\Spc}(K(\Pi,1))$, and therefore the two 
$B$--structures $\nu$ and $\nu'$ are cobordant.
Theorem 2 in~\cite{K} now implies the existence of numbers $t,s \geq 0$
and of a diffeomorphism
\[
\Theta \co X \#   t(S^2 \times S^2) \longrightarrow X'   \# s(S^2 \times S^2)
\]
compatible with the normal $B$--structures obtained by $\nu$ respectively $\nu'$ and the
canonical $B$--structure on $S^2 \times S^2$.
In particular, $\Theta$ maps $c_1(s')$ to $c_1(s)=PD(\xi)$. By construction, $c_1(s')$ has support in the
simply connected part $Z$ and can therefore be represented by a $\pi_1$--null embedding. Pulling back this
surface via $\Theta$ we obtain a $\pi_1$--null embedding $F \rightarrow X \# t(S^2 \times S^2)$ representing $\xi$
as desired.

As to the general case, assume that $\xi \in H_2(X;\Zset)$ is a spherical and characteristic homology class.
Denote the usual generator of $H_2(\CP^2;\Zset)$ by $\gamma$ and consider the 4--manifold
$X'=X \# \CP^2$. Let $\xi'=(\xi,\gamma) \in H_2(X';\Zset)$. The homology class of a generically
embedded $\CP^1$ in $\CP^2$ is spherical and has intersection number one with $\xi'$. By what we just proved, this implies
that $\xi'$ can be stably represented by a $\pi_1$--null embedding. As $\gamma$ is characteristic, Proposition~\ref{embeddings}
can be applied and we obtain that also $\xi$ can be stably represented by a $\pi_1$--null embedding.
\end{proof}

\begin{proof}[Proof of Theorem~\ref{FKlike}]\label{proof}
By Theorem~\ref{kem}, a characteristic class $\xi$ which can be stably
represented by an embedded sphere fulfills $\xi \cdot \xi \equiv \sigma(X) \mod 16$.
Let us now assume that conversely, the
congruence  $\xi \cdot \xi \equiv \sign(X) \mod 16$ holds.
By Theorem~\ref{charcase} 
we can assume that the class $\xi$ can be represented
by a smoothly embedded surface $F \subset X$ such that $\pi_1(F) \rightarrow \pi_1(X)$
is trivial.

By Theorem 1 in \cite{FK}, the Arf invariant of $F$ vanishes. 
Now we could proceed  as in the proof of Theorem
2 in \cite{FK} to obtain an embedded sphere
in \mbox{$X \# m(S^2 \times S^2)$}
for some $m$ representing $[F]$. However,  there is a slightly different argument using 
transversal spheres instead of framed surgery along circles.

Since the Arf invariant of $F$ is zero, we can find a homologically non--trivial circle 
$C \subset F$ with $q(C)=0$, here our notation is the same as in~\cite{FK}.
Pick an embedded disk $D \subset X$ with boundary $C$ whose interior intersects the surface $F$ transversely with 
algebraic intersection number $k$.
After performing boundary twists, we can assume that $k=0$.
Fix a framing of $C$
in $F$ and let $l$ denote the integer valued obstruction to extending this framing to a section of the
normal bundle of $D$. According to the definition of $q$, we have $k + l \equiv q(C)=0 \mod 2$, hence
$l$ is even.

Now we will use transversal spheres 
to remove the intersection points of $F \setminus C$ and $D$, as follows.  After passing
from $X$ to $X \# (S^2 \times S^2)$ and attaching  one of the factors to $D$
(note that this does note change $l$, since the attached sphere has trivial normal bundle), 
we can assume that $D$ 
has a transversal sphere which does not meet $F$, 
namely the second factor of $S^2 \times S^2$.
We can use parallel copies of this transversal sphere to remove all the intersection points between the interior of 
$D$ and $F$. As $k=0$, the homology classes of all these copies add up to zero.
Hence we obtain a new surface $F'$ which still
has homology class $[F]$ and genus $g(F)$ such that $D$ intersects $F'$ only in $C$. Since the surface $F'$
is obtained from $F$ by cutting out disks and gluing in other disks instead, the circle $C$ 
is still non--trivial.
Furthermore, we did not change the disk $D$, and hence the obstruction to extending a framing of $C$
over $D$ is still~$l$.

The arguments given so far show that we can arrange for the interior of the disk $D$ to be disjoint from $F$.
After adding another copy of $S^2 \times S^2$, we can find an embedded sphere $S$ with self
intersection number $-l$ disjoint from $F$ and $D$, note that $l$ is even. Tubing this sphere 
into the disk $D$ gives a disk $D'$ with boundary $C$, still disjoint from $F$, such that
the framing obstruction  vanishes.
By doing surgery along $C$, we can now obtain a new surface representing $[F]$
with genus $g(F)-1$. Repeating the argument, we can therefore construct an embedded sphere as desired. 
 \end{proof}

\section{The topological case}\label{top}

Up to now, we have always been working in the smooth category, i.e.\ 
we have considered smooth 4--manifolds, smoothly
immersed spheres and smooth embeddings. However, it turns out that most of our results ---
including Theorem~\ref{charcase} --- are also true in the topological category.

First of all, note that --- as explained in \cite{FQ} --- every homotopy class of maps from a surface
to a topological 4--manifold contains an immersion. We can define self intersection numbers for these
immersions as in the smooth case, and again obtain a map
\[
\bar{\mu} \co \pi_2(X) \longrightarrow \Zset[\pi_1(X)] / ( \langle \alpha - \alpha^{-1} \rangle + \Zset )
\]
for every topological 4--manifold $X$.
Lemma~\ref{pi1null} then has an obvious generalisation to the topological case, and the definition of the
stable self intersection numbers still makes sense. Since the proof of Proposition~\ref{embeddings}
is based on algebraic arguments, it can easily be adapted to the topological case to obtain

\begin{prop}\label{topembeddings}
Suppose we are given a topological 4--manifold $X$ and a homology 
class $\xi \in H_2(X;\Zset)$.
Then the following conditions are equivalent.
\begin{enumerate}
\item There is an immersion $f \co S^2 \rightarrow X$ representing $\xi$ 
such that $\mu_s(f)=0$.
\item There is, for some $k$, a locally flat $\pi_1$--null embedding 
$F \rightarrow X \# k(S^2 \times S^2)$
representing the homology class $(\xi,0, \dots, 0)$.
\item There is, for some simply connected topological 4--manifold $Y$, a locally flat
$\pi_1$--null embedding $F \rightarrow
X \# Y$ such that the homology class of $F$ is $(\xi,c)$ with a characteristic class $c \in H_2(Y;\Zset)$.
\end{enumerate}
\end{prop}

Now it turns out that also Theorem~\ref{charcase} remains true in the topological category. Instead of going
through the proof and checking carefully at which points the smooth structure is used, it is
more convenient to reduce the topological case to the smooth case using an additional argument.

\begin{thm}\label{topcharcase}
A characteristic class in a topological 4--manifold can be stably represented by
a locally flat $\pi_1$--null embedding if and only if it is spherical.
\end{thm}

\begin{proof} 
Again it is clear that a class which can be stably represented by a $\pi_1$--null embedding
is spherical. Now let $X$
be a 4--manifold with a characteristic spherical class $\xi \in H_2(X;\Zset)$.
First let us suppose that the Kirby--Siebenmann invariant $\ks(X)$ is zero. Then there is some $t$ such
that the manifold $X'=X \# t(S^2 \times S^2)$ is smoothable. Applying Theorem~\ref{charcase} to this
manifold shows that the class $(\xi,0) \in H_2(X';\Zset)$ can be stably represented by a $\pi_1$--null
embedding. By the very definition of ``stably representable'', the same is true for $\xi$.

In the case that $\ks(X)=1$, consider the 4--manifold $X'=X \# E_8$ and the characteristic class
$\xi'=(\xi,0)$. Then $\ks(X')=\ks(X)+1=0$, and by the first part of the proof, 
we can represent the class $\xi'$ stably
by a locally flat $\pi_1$--null embedding. By Proposition~\ref{topembeddings}, the same is true for $\xi$. 
\end{proof}

Finally let us prove an analogue of Theorem~\ref{FKlike} in the topological case.
If $F \subset X$ is a topologically locally flat embedded surface in a 4--manifold $X$
such that $\pi_1(F) \rightarrow \pi_1(X)$ is trivial, one can define a quadratic form
\[
q \co H_1(F;\Zset_2) \longrightarrow \Zset_2
\] 
as in \cite{FK}. However, it is no longer true that
$\Arf(q)=\frac{1}{8}(F \cdot F - \sign(X))$ as the example of a trivially embedded sphere in
$E_8$ shows. Instead, we have 
\[
\Arf(q) = \frac{1}{8}(F \cdot F - \sign(X)) + \ks(X) \mod 2,
\]
where $\ks(X)$ is the Kirby--Siebenmann smoothing obstruction, see \cite{LW}
(there this is only proved in the case that $X$ is simply connected, but one can always
reduce to this case by surgery along embedded circles).
In particular, we obtain that 
$
\frac{1}{8}(F \cdot F - \sign(X)) \equiv \ks(X) \mod 2
$
if the surface $F$ is a sphere. If $\Arf(F)=0$, then we can obtain an embedded
sphere in $X \# k(S^2 \times S^2)$ for some $k$ representing $[F]$ as in the proof of Theorem~\ref{FKlike}. 
Hence Theorem~\ref{topcharcase} implies the following.

\begin{thm}
Let $X$ be a topological 4--manifold  and let $\xi \in H_2(X;\Zset)$ be a spherical and
characteristic homology class. Then $\xi$ can be stably represented by a locally flat embedding of a  2--sphere
if and only if 
\[
\frac{1}{8}(\xi \cdot \xi - \sign(X)) \equiv \ks(X) \mod 2,
\]
where $\ks(X) \in \Zset_2$ denotes
the Kirby--Siebenmann smoothing obstruction.
\end{thm}

This result was obtained in~\cite{LW} in the special case that the manifold $X$ is simply connected,
note that in this case, every homology class is spherical.

\Addressesr
\end{document}
